\documentclass{amsart}
\usepackage{amssymb}
\usepackage{amscd}
\usepackage{latexsym}
\usepackage{amsfonts}
\usepackage{amsmath}
\newtheorem{theorem}{Theorem}
\newcommand{\eproof}{\begin{flushright} $\square$ \end{flushright}}

\newcommand{\tr}{\mathop{\fam0 Tr}\nolimits}

\newcommand{\End}{\mathop{\fam0 End}\nolimits}

\newcommand{\Aut}{\mathop{\fam0 Aut}\nolimits}

\newcommand{\Id}{\mathop{\fam0 Id}\nolimits}

\newcommand{\bC}{{\mathbb C}}
\newcommand{\bR}{{\mathbb R}}
\newcommand{\C}{C}

\newcommand{\Z}{{\mathbb Z}}

\newcommand{\ra}{\mathop{\fam0 \rightarrow}\nolimits}

\newcommand{\cH}{ {\mathcal H}}
\newcommand{\T}{ {\mathcal T}}

\newcommand{\V}{ {\mathcal H}}
\renewcommand{\L}{{\mathcal L}}
\newcommand{\Sk}{{\mathcal S}}

\newcommand{\cZ}{ {\mathcal Z}}

\renewcommand{\P}{ {\mathbb P}}

\newcommand{\s}{\sigma}

\newcommand{\Teim}{Teichm{\"u}ller }

\begin{document}
\title[Asymptotics of the Hilbert-Smith norm in TQFT]
{Asymptotics of the Hilbert-Smith norm of\\ curve operators in
TQFT}
\author{J\o rgen Ellegaard Andersen}
\address{Department of Mathematics\\
        University of Aarhus\\
        DK-8000, Denmark}
\email{andersen@imf.au.dk}

\begin{abstract}
Applying standard techniques from Toeplitz operator theory, we
analyze the asymptotics of the Hilbert-Smith norms of the TQFT
operators coming from isotopy classes of one dimensional
oriented submanifolds on a closed oriented surface. We thereby obtain a
Toeplitz operator interpretation and generalization of the asymptotic formula
obtained recently by March\'{e} and Narimannejad in \cite{MN}.
\end{abstract}

\maketitle

\section{Introduction}

By the very axioms of a general TQFT (see \cite{T}), we have
the following setup.

Any TQFT $\cZ$ associates a finite dimensional Hermitian vector
space $\cZ(\Sigma)$ to a closed oriented surface $\Sigma$. It also
induces a projective representation
$$ \cZ : \Gamma \ra \Aut(\P(\cZ(\Sigma))),
$$
where $\Gamma$ is the mapping class group of $\Sigma$.

For any tripel $(M,L,\lambda)$, where $M$ is compact oriented three
manifold,  $L$ is a one dimensional oriented framed sub-manifold
in $M-\partial M$ and $\lambda$ is a labeling of each of the
components of $L$ by an element from a certain finite label set, which is
part of the theory $\cZ$, one get a vector
$$
\cZ(M,L,\lambda) \in \cZ(\partial M).
$$

What interest us here is the operators associated to labeled simple
closed curves. More generally suppose $\gamma$ is an oriented one
dimensional sub-manifold $\gamma$ of $\Sigma$ and that $\lambda$ is
a labeling of each of the components of $\gamma$. Then by embedding
$\gamma$ in $\Sigma\times [0,1]$ via the identification of $\Sigma$
with $\Sigma\times\{1/2\}$ and using the "black board" framing along $\Sigma\times\{1/2\}$,
we get the {\em curve operator}
$$
\cZ(\gamma,\lambda) \in \End(\cZ(\Sigma)).
$$

Using the above mentioned Hermitian structure we have the
Hilbert-Smith norm or inner product in its polarized form:
$$\langle\cZ(\gamma_1,\lambda_1), \cZ(\gamma_2,\lambda_2)\rangle =
\tr(\cZ(\gamma_1,\lambda_1)\cZ(\gamma_2,\lambda_2)^*).$$

In this paper we will compute the large level asymptotics of this
Hilbert-Smith Norm for the Reshetikhin-Turaev quantum $SU(n)$ TQFT
$Z^{(k)}$ with the use of Toeplitz operator theory. Let us make this precise.

The TQFT's $Z^{(k)}$ was rigorously constructed by Reshetikhin and
Turaev in \cite{RT1} and \cite{RT2} using the theory of quantum
groups. A purely skein theory model for these TQFT's was developed
in \cite{BHMV1}, \cite{BHMV2} and \cite{B1}.

The label set $\Lambda_k$ of the theory $Z^{(k)}$ is given as follows:
Let $\Lambda$ be the set of finite dimensional irreducible
representations of $SU(n)$. By the theory of dominant weight modules,
$\Lambda$ is identified with the set of dominant
weights of the Lie algebra of $SU(n)$. The label set for $Z^{(k)}$ is the
following subset of $\Lambda$:
$$
\Lambda_k = \{ \lambda \in \Lambda  |  ( \lambda, \theta ) \leq k\},
$$
where $( \cdot, \cdot )$ is the Cartan-Killing form and
$\theta$ is the highest root of the
Lie algebra of $SU(n)$ (and $( \cdot, \cdot )$ is normalized by requiring
that $( \theta, \theta ) = 2$).
We observe that any $\lambda\in \Lambda$ is contained in $\Lambda_k$
for $k$ large enough. The operation of taking the dual
representation gives an involution of $\Lambda$, which preserves
$\Lambda_k$ for each $k$.

Let $n\geq 2$ and consider the space $Z^{(k)}(\Sigma)$ associated
to a surface of genus $g>1$. In fact, we will consider a slightly
more general situation, namely where we label a point $x$ on
$\Sigma$ by an element from $\Lambda_k$. We will only
consider the following situation: Let $\lambda_0\in \Lambda_1$ be
 the unique element with the
property that the exponential map applied to
its associated co-weight generate the center of $SU(n)$. Then
let $d \in\Z/n\Z$, label $x$ with $d\lambda_0$ and let
$Z^{(k)}(\Sigma)$ be the vector space $Z^{(k)}$ associates to
$\Sigma$ label by $d\Lambda_0$ at $x$. In case $d=0$, this is
simply just the vector space $Z^{(k)}$ associates to the surface $\Sigma$
without any labeled points.

Suppose now $\gamma$ is a one dimensional submanifold of $\Sigma$ and
that $\lambda$ is a labeling of the components
of $\gamma$ by finite dimensional representations.
For large enough $k$, we therefore get the curve operator
$$
Z^{(k)}(\gamma, \lambda) \in \End(Z^{(k)}(\Sigma)).
$$
If all components of $\gamma$ are labeled by the defining
representation of $SU(n)$, we just write $Z^{(k)}(\gamma)$ for the corresponding
curve operator. In section \ref{sec1} we review the construction of the curve operators
$Z^{(k)}(\gamma,\lambda)$.

We now proceed to discuss the semiclassical limit of these sequences
of curve operators, namely holonomy functions on the $SU(n)$ moduli
space.

 Let $M$ be the moduli space of flat
$SU(n)$-connections on $\Sigma - x$ with holonomy $d\in \Z/n\Z \cong
Z_{SU(n)}$ around $x$.

Suppose now that each of the components of $\gamma$ is labeled by a
finite dimensional representation of $SU(n)$ as discussed above.
Then we have the holonomy function $h_{\gamma,\lambda}$ associated to $\gamma$
defined on
$M$ by taking the product over the
components of $\gamma$ of the trace in the corresponding
representation of the holonomy around the component. Note that
$h_{\gamma,\lambda}$ only depends on the free homotopy class of
$\gamma$. Further $h_{\gamma,\lambda}$ is constant if $\gamma$ is
nul-homotopic. If we label all components of $\gamma$ by the
defining representation of $SU(n)$, we simply denote the corresponding
function by $h_\gamma$.

On the holonomy functions we have the $L_2$-inner product:
$$\langle h_{\gamma_1,\lambda_1}, h_{\gamma_2,\lambda_2}\rangle = \frac1{m!}\int_M
h_{\gamma_1,\lambda_1}{\overline h_{\gamma_2,\lambda_2}} \omega^m,$$
where $m = (g-1)(n^2 -1)$ and $\omega$ is the symplectic form on $M$.

The main purpose of this paper is to explain how the following theorem
can be deduced from standard results from the theory of Toeplitz operators and by
the approximation result Theorem \ref{MainA6} proved in \cite{A6}.

\begin{theorem}\label{Main}
For all pairs of oriented one dimensional submanifolds $\gamma_1, \gamma_2$ of $\Sigma$ and
all labelings $\lambda_1, \lambda_2$ of their components by finite
dimensional irreducible representations of $SU(n)$, we have that
$$
\langle h_{\gamma_1,\lambda_1}, h_{\gamma_2,\lambda_2}\rangle =
\lim_{k\ra \infty} k^{-(g-1)(n^2-1)} \langle Z^{(k)}(\gamma_1,
\lambda_1),Z^{(k)}(\gamma_2,\lambda_2)\rangle,$$

\end{theorem}

To explain how this statement is related to Toeplitz operators, let us
review the $SU(n)$ gauge theory construction of $Z^{(k)}(\Sigma)$.

By applying geometric quantization at level $k$ to the moduli space
$M$ one gets a vector bundle $\V^{(k)}$ over \Teim space $\T$. The
fiber of this bundle over a point $\sigma\in \T$ is $\V^{(k)}_\sigma
= H^0(M_{\sigma},\L_\s^k)$, where $M_{\s}$ is $M$ equipped with a
complex structure induced from $\sigma$ and $\L_\s$ is an ample
generator of the Picard group of $M_{\s}$.

The main result pertaining to this bundle $\V^{(k)}$ is that its
projectivization $\P(\V^{(k)})$ supports a natural flat connection.
This is a result proved independently by Axelrod, Della Pietra and
Witten \cite{ADW} and by Hitchin \cite{H}. Now, since there is an
action of the mapping class group $\Gamma$ of $\Sigma$ on $\V^{(k)}$
covering its action on $\T$, which preserves the flat connection in
$\P(\V^{(k)})$, we get for each $k$ a finite dimensional projective
representation, say $Z^{(k)}$, of $\Gamma$, namely on the covariant
constant sections, $\P(Z^{(k)}(\Sigma))$ of $\P(\V^{(k)})$ over
$\T$. This sequence of projective representations $Z^{(k)}$, $k\in
{\mathbb N}$ is the {\em quantum $SU(n)$ representations} of the
mapping class group $\Gamma$.

For each $f\in C^\infty(M)$ and each point $\sigma\in \T$, we have
the {\em Toeplitz operator}
\[T^{(k)}_{f,\sigma} : H^0(M_{\sigma},\L_\s^k) \ra H^0(M_{\sigma},\L_\s^k)\]
which is given by
\[T^{(k)}_{f,\sigma} = \pi^{(k)}_\sigma(fs)\]
for all $s\in H^0(M_{\sigma},\L_\s^k)$, where $\pi^{(k)}_\sigma$ is
the orthogonal projection onto $H^0(M_{\sigma},\L_\s^k)$. We get
smooth section of $\End(\V^{(k)})$ over $\T$
\[T_{f}^{(k)} \in C^\infty(\T,\End(\V^{(k)})) \]
by letting $T_{f}^{(k)}(\sigma) = T_{f,\sigma}^{(k)}$. See also
\cite{A3}  and \cite{A4} where we use these operators to prove the asymptotic faithfulness of the
quantum $SU(n)$ representations.

The $L_2$-inner product on $\C^\infty(M,\L^k)$ induces an inner
product on $H^0(M_\sigma, \L_\sigma^k)$, which in turn induces the
operator norm $\|\cdot\|$ on $\End(H^0(M_\sigma, \L_\sigma^k))$.
Hence for any $A\in C^\infty(\T,\End(\V^{(k)}))$ we get the smooth
function $\|A\|$ on $\T$.

\begin{theorem}\label{MainToeplitz}
For any two smooth functions $f,g \in \C^\infty(M)$ one has that
$$\langle f, g\rangle = \lim_{k\ra \infty} k^{-m} \tr(T^{(k)}_{f} (T^{(k)}_{g})^*),$$
where the real dimension of $M$ is $2m$.
\end{theorem}

We give a proof of this Theorem in section \ref{sec2}. It is a well
known result in the field of Toeplitz operators.

As we will see in section \ref{sec2}, Theorem \ref{Main} follows from Theorem \ref{MainToeplitz} and the
following Theorem:

\begin{theorem}\label{MainA6}
For any one dimensional oriented submanifold $\gamma$ and any
labeling $\lambda$ of the components of $\gamma$, we have that
$$
\lim_{k \ra \infty} \|Z^{(k)}(\gamma, \lambda) -
T^{(k)}_{h_{\gamma,\lambda}}\| = 0.
$$
\end{theorem}

The proof of this theorem, which we give in \cite{A6}, uses only gauge theory and then
asymptotic analysis of certain leafwise elliptic operators on the
moduli space, plus ideas from conformal field theory. However, it
does {\em not} rely on the our joint work with K. Ueno \cite{AU1},
\cite{AU2} and \cite{AU3}.

Let us now specialize to the setting discussed by March\'{e} and
Narimannejad. The TQFT $V_p$ is the one defined by Blanchet,
Habegger, Masbaum and Vogel in \cite{BHMV1} and \cite{BHMV2}, where
$p = 2 r$, $r$ being an integer and we choose the $2p$'th root of $1$
to be $A = - e^{i \pi/2r}$. I.e. this is the Skein theory model for the quantum
$SU(2)$ Reshetikhin-Turaev TQFT.

There is also an inner product on $V_p(\Sigma)$ and the
representation $V_p$ of $\Gamma$ on this vector space is unitary.

Let $\Sk(\Sigma)$ be the free $\bC$-vector space generated by
isotopy classes of one-dimensional sub-manifolds of $\Sigma$. We do
not need any orientation in this case since this is the $SU(2)$
theory. The vector space $\Sk(\Sigma)$ is known to be isomorphic to the
vector space of algebraic functions on the $SL(2,\bC)$ moduli space by the work of
Bullock \cite{Bul} and the work of Przytycki and Sikora \cite{PS}.

As above we get a sequence of linear maps
\[V_p : \Sk(\Sigma) \ra \End(V_p(\Sigma)),\]
and the inner product on $V_p(\Sigma)$ induces an inner product
$\langle \cdot,\cdot \rangle_p$ on $\End(V_p(\Sigma))$.

In \cite{MN} March\'{e} and Narimannejad proved

\begin{theorem}[March\'{e} and Narimannejad]\label{MNT}
For all $\gamma_1,\gamma_2\in \Sk(\Sigma)$
$$
\langle h_{\gamma_1}, h_{\gamma_2}\rangle = \lim_{p\ra \infty} 2^{(3g - 3)}p^{-(3g - 3)} \langle V_p(\gamma_1),V_p(\gamma_2) \rangle_p ,$$
where $M$ is the moduli space of flat $SU(2)$-connections, i.e. $(n,d)= (2,0)$.
\end{theorem}

We have in this paper interpreted this formula in terms of Toeplitz
operators and in the process generalize the formula to the
$SU(n)$-case. In fact we have given an alternative proof of their theorem, since
we can deduce it from Theorem \ref{Main}:

Our joint work with K. Ueno \cite{AU1}, \cite{AU2} and \cite{AU3},
combined with the work of Laszlo \cite{La1} gives the following
result for the $(n,d) = (2,0)$ theory:
\begin{theorem}[AU]\label{projiso}
There is a projective linear isomorphism of representations
of $\Gamma$
\[I_k :  \P(V_{2k+4}(\Sigma)) \ra \P(Z^{(k)}(\Sigma)).\]
\end{theorem}

The projective linear isomorphism $I_k$ induces a  linear algebra
isomorphism of representations of $\Gamma$
\[I^e_k :  \End(V_{2k+4}(\Sigma)) \ra  \End(Z^{(k)}(\Sigma)).\]

From the construction of $I_k$ it follows that
\begin{theorem} \label{isocurveop}
For any one dimensional submanifold $\gamma$ of $\Sigma$ we have
that
\[I_k^e(V_{2k+4}(\gamma))= Z^{(k)}(\gamma).\]
\end{theorem}

We see that Theorem \ref{isocurveop} and Theorem \ref{Main}
implies March\'{e} and Narimannejad's Theorem \ref{MNT}: Using the
the algebra isomorphism $I^e_k$, we compute that
\begin{eqnarray*}
\tr(V_{2k +4}(\gamma_1) V_{2k+4}(\gamma_2)) & = & \tr(I^e_k(V_{2k +4}(\gamma_1) V_{2k +4}(\gamma_2)))\\
& =&  \tr(I^e_k(V_{2k +4}(\gamma_1)) I^e_k(V_{2k +4}(\gamma_2)))\\
& = & \tr(Z^{(k)}(\gamma_1)Z^{(k)}(\gamma_2)).
\end{eqnarray*}
Now we apply Theorem \ref{Main} to get Theorem \ref{MNT} as a special case,
since the curve operators $V_{2k+4}(\gamma)$ are self-adjoint.

This paper is organized as follows. In section \ref{sec1} we review the construction of the curve operators $Z(\gamma, \lambda)$,  and
in section \ref{sec2} we prove
theorem \ref{Main}.

We would like to thank Gregor Masbaum for discussions related to this work.

\section{Construction of the curve operators}\label{sec1}

In this section we review how the curve operators $Z^{(k)}(\gamma,\lambda)$
are constructed for the theory $Z^{(k)}$.

Let $\Sigma'$ be the surface with boundary obtained from cutting
$\Sigma$ along $\gamma$. Then any labeling $\mu$ of $\gamma$ by elements
from $\Lambda_k$ induces a labeling of the boundary components of
$\Sigma'$, using the following convention: If $\Sigma'$ induces
the same orientation on a component of $\partial\Sigma'$ as
$\gamma$ does, we use the label for that component, else we use the
dual of the label. We also denote this labeling of $\partial\Sigma'$ by
$\mu$.

Since $Z^{(k)}$ is also a modular functor one
can factor the space $Z^{(k)}(\Sigma)$ into a direct sum "along"
$\gamma$ as a sum over all labelings of $\gamma$. That is we get an isomorphism

\begin{equation}
Z^{(k)}(\Sigma) \cong \bigoplus_{\mu} Z^{(k)}(\Sigma',\mu).\label{Fact}
\end{equation}
The sum here runs over all labelings $\mu$ of $\gamma$ by elements from
$\Lambda_k$. Strictly speaking we need base points on all the components of
$\partial \Sigma'$ to define $Z^{(k)}(\Sigma',\mu)$. However, the corresponding
subspaces of $Z^{(k)}(\Sigma)$ does not depend on the choice of base points.
 The isomorphism (\ref{Fact}) induces an isomorphism
\begin{equation*}
\End(Z^{(k)}(\Sigma)) \cong \bigoplus_{\mu}
\End(Z^{(k)}(\Sigma',\mu)).
\end{equation*}
which also induces a direct sum decomposition of $\End(Z^{(k)}(\Sigma))$ which
is independent of the choice of the base points.

The curve operator $Z^{(k)}(\gamma,\lambda)$ is diagonal with respect
to this direct sum decomposition along $\gamma$. One has the
formula
\[Z^{(k)}(\gamma,\lambda) = \bigoplus_{\mu} S_{\lambda,\mu}
(S_{0,\mu})^{-1}\Id_{Z^{(k)}(\Sigma',\mu)}.\]
Here
\[S_{\lambda,\mu} = \prod_{i=1}^s S_{\lambda(\gamma_i),\mu(\gamma_i)}\]
where $\gamma_i$, $i=1,\ldots, s$ are the components of $\gamma$,
and $S_{\lambda_1,\lambda_2}$, $\lambda_1,\lambda_2 \in \Lambda_k$
is the $S$-matrix of the theory $Z^{(k)}$.
See e.g. \cite{B1} for a derivation of this.

In the gauge theory picture, the decomposition (\ref{Fact}) is obtained
as follows (see \cite{A6} for the details):

One considers a one parameter family of complex structures $\sigma_t\in
\T$, $t\in \bR_+$, such that the corresponding family in the moduli space of
curves converges in the Mumford-Deligne boundary to a nodal curve,
which topologically corresponds to shrinking $\gamma$ to a point.
By the results of \cite{A1.5} the corresponding sequence of
complex structures on the moduli space $M$ converges to a non-negative
polarization on $M$, whose isotropic foliation is spanned by the Hamiltonian
vector fields associated to the holonomy functions of $\gamma$.
The main result of \cite{A6} is that the covariant constant sections
of $\cH^{(k)}$ along the family $\sigma_t$ converges to
distributions supported on the Bohr-Sommerfeld leaves
of the limiting non-negative polarization as $t$ goes to infinity.
The direct sum of the geometric quantization of the level $k$
Bohr-Sommerfeld levels of the limit non-negative polarization is
precisely the left hand side of (\ref{Fact}). A sewing-construction
inspired from conformal field theory (see \cite{TUY}) is then applied to show that
the resulting linear map from the right hand side of (\ref{Fact})
to the left hand side is an isomorphism. This is described in
detail in \cite{A6}.

Theorem \ref{MainA6} follows from this asymptotic
study: One considers the explicit expression for the $S$-matrix,
as given in formula (13.8.9) in Kac's book \cite{Kac}
\begin{equation}
S_{\lambda,\mu}/S_{0,\mu} =  \lambda(e^{-2 \pi i \frac{\check{\mu} +
\check{\rho}}{k+n}}),\label{KacS}
\end{equation}
where $\rho$ is half of the sum of the positive roots and
$\check{\nu}$ ($\nu$ any element of $\Lambda$) is the unique element of the Cartan subalgebra of the Lie algebra of
$SU(n)$ which is dual to $\nu$ with respect to the Cartan-Killing form $(\cdot,\cdot)$.

From the expression (\ref{KacS}) one sees that under the isomorphism
$\check{\mu}\mapsto \mu$,
$S_{\lambda,\mu}/S_{0,\mu}$  makes sense for
any $\check{\mu}$ in the Cartan subalgebra of the Lie algebra of
$SU(n)$. Furthermore one finds that the values of this sequence of functions
(depending on $k$) is asymptotic to the values of the holonomy function
$h_{\gamma,\lambda}$ at the level $k$ Bohr-Sommerfeld sets of the limiting
non-negative polarizations discussed above (see \cite{A1}). From this one can deduce Theorem \ref{MainA6}.
See again \cite{A6} for details.

\section{Toeplitz operators and proof of Theorem \ref{Main}.}\label{sec2}

Let us recall basics on the theory of Toeplitz operators and apply
this theory to prove Theorem \ref{Main}.

On $\C^\infty(M,\L^k)$ we have the $L_2$-inner product:
\[\langle s_1, s_2 \rangle = \frac{1}{m!}\int_M (s_1,s_2) \omega^m\]
where $s_1, s_2 \in \C^\infty(M,\L^k)$ and $(\cdot,\cdot)$ is the
fiberwise Hermitian structure in $\L^k$.

Now let $\sigma\in \T$. Then this $L_2$-inner product gives the
 orthogonal projection $$\pi^{(k)}_\sigma : \C^\infty(M,\L^k) \ra
 H^0(M_\sigma,\L_\sigma^k).$$
For each $f\in \C^\infty(M)$ consider the associated {\em Toeplitz
operator} $T_{f,\sigma}^{(k)}$ given as the composition of the
multiplication operator $$M_f : H^0(M_\sigma,\L_\sigma^k) \ra
\C^\infty(M,\L^k)$$ with the orthogonal projection:
\[T_{f,\sigma}^{(k)}(s) = \pi^{(k)}_\sigma(f s).\]
Then $T_{f,\sigma}^{(k)} \in \End(H^0(M_\sigma,\L_\sigma^k))$, and
we get a smooth section
\[T_{f}^{(k)} \in C^\infty(\T,\End(\V^{(k)})) \]
by letting $T_{f}^{(k)}(\sigma) = T_{f,\sigma}^{(k)}$ (see
\cite{A3}).

The $L_2$-inner product on $\C^\infty(M,\L^k)$ induces an inner
product on $H^0(M_\sigma, \L_\sigma^k)$, which in turn induces the
operator norm $\|\cdot\|$ on $\End(H^0(M_\sigma, \L_\sigma^k))$.

We need the following Theorem on Toeplitz operators
due to Bordemann, Meinrenken and Schlichenmaier (see \cite{BMS}, \cite{Sch},
\cite{Sch1} and \cite{Sch2}).
\begin{theorem}[Bordemann, Meinrenken and Schlichenmaier]\label{BMS1}
For any $f\in \C^\infty(M)$ we have that
\[\lim_{k\ra \infty}\|T_{f,\sigma}^{(k)}\| = \sup_{x\in M}|f(x)|.\]
\end{theorem}
Since the association of the sequence of Toeplitz operators
$T^k_{f,\sigma}$, $k\in \Z_+$ is linear in $f$, we see from this
Theorem, that this association is faithful.

We also have the following two theorems from \cite{BMS}:

\begin{theorem}[Bordemann, Meinrenken and Schlichenmaier]\label{BMS2}
For any $f, g\in \C^\infty(M)$ we have that
\[\|T_{f,\sigma}^{(k)}T_{g,\sigma}^{(k)} - T_{fg,\sigma}^{(k)}\| = O(k^{-1}).\]
\end{theorem}

\begin{theorem}[Bordemann, Meinrenken and Schlichenmaier]\label{BMS3}
For any $f\in \C^\infty(M)$ we have that
\[\lim_{k\ra \infty}k ^{-m} \tr(T_{f,\sigma}^{(k)}) = \frac1{m!} \int_M f \omega^m.\]
\end{theorem}

\proof[Proof of Theorem \ref{MainToeplitz}]
It is immediate from the definition of a Toeplitz operator that
$$(T^{(k)}_{f,\sigma})^* = T^{(k)}_{{\overline f},\sigma}.$$
But then by Theorem \ref{BMS3} and \ref{BMS2}
\begin{eqnarray*}
\langle f,g \rangle & = &\lim_{k\ra \infty}k ^{-m} \tr(T_{f{\overline g},\sigma}^{(k)})\\
& = &\lim_{k\ra \infty}k ^{-m} \tr(T_{f,\sigma}^{(k)}T_{{\overline g},\sigma}^{(k)})\\
& = &\lim_{k\ra \infty}k ^{-m} \tr(T_{f,\sigma}^{(k)}(T_{g,\sigma}^{(k)})^*)).
\end{eqnarray*}
\eproof

Using Theorem \ref{MainA6} we can now give a Toeplitz operator proof of Theorem \ref{Main}.

\proof[Proof of Theorem \ref{Main}]
By Theorem \ref{MainA6} we get that
$$\lim_{k\ra \infty}\|T_{h_{\gamma_1,\lambda_1}}^{(k)}T_{h_{\gamma_2,\lambda_2}}^{(k)} - Z^{(k)}(\gamma_1,\lambda_1) Z^{(k)}(\gamma_2,\lambda_2)\| = 0.$$
But now recall that the Hilbert-Smith norm is bounded by the square root of the dimension of $H^0(M_\sigma, \L_\sigma^k)$ times the operator norm.
Since the dimension of $H^0(M_\sigma, \L_\sigma^k)$ grows like $k^m$, we certainly get that
$$\lim_{k\ra \infty}k ^{-m}|\tr(T_{h_{\gamma_1}}^{(k)}(T_{h_{\gamma_2}}^{(k)})^*) -
\tr(Z^{(k)}(\gamma_1,\lambda_1) (Z^{(k)}(\gamma_2,\lambda_2))^*)| =0.$$
But now the desired result follows from Theorem
\ref{MainToeplitz}.
\eproof

\end{document}